\documentclass[12pt]{amsart}
\usepackage[T2A]{fontenc}
\usepackage[cp1251]{inputenc}
\usepackage[english]{babel}
\usepackage{amsmath,latexsym,amsthm,amsfonts,tipa,upgreek}
\usepackage{amssymb,amscd,graphpap, stmaryrd}

\begin{document}

\title{Classification of  spaces of Continuous Function on ordinals}
\author{Genze L.V., Gul'ko S.P., Khmyleva T.E.}

\thanks{The study was carried out with the financial support of the Russian Foundation for Basic
Research in the framework of the scientic project N 17-51-18051 }

\maketitle

\begin{abstract}
We conclude the classification of spaces of continuous functions on ordinals carried out by Gorak \cite{Gorak}.
This gives a complete topological classification of the spaces $C_p([0,\alpha])$ of all continuous real-valued functions on compact segments of ordinals endowed with the topology of pointwise convergence. 
Moreover, this topological classification of the spaces  $C_p([0,\alpha])$ completely coincides with their uniform classification.
\vspace{3 mm}

\vspace{3 mm}


\end{abstract}

\begin{center}
\textbf{1. Introduction}
\end{center}

Our terminology basically follows \cite{Engel}. In particular, we understand cardinals as initial ordinals, compare \cite{Engel}, page 6.
A segment of the ordinals $[0, \alpha]$ is endowed with a standard order topology.
The symbol $C_p ([0, \alpha])$ denotes the set of all continuous real-valued functions defined on $[0, \alpha]$ and endowed with the topology of pointwise convergence.

A complete linear topological classification of Banach spaces $C([0,\alpha])$ was carried out in \cite{GulOs} and independently in \cite{Kisl} (for the initial part of this classification, see also \cite{BesPel} and \cite{Sem}).
Similar complete linear topological classification for $C_p([0, \alpha])$ can be found in \cite{Gulko, BaarsDeGroot}.

The topological classification of the spaces $C_p ([0, \alpha])$ is carried out in the Gorak's paper \cite{Gorak}, in which the question whether the spaces $C_p ([0,\alpha])$ and $ C_p ([0, \beta]) $ are homeomorphic is solved for all ordinals $\alpha$ and $\beta$ with except for the case $\alpha = k ^+ \cdot k $, $\beta = k ^+ \cdot k^+$, where $ k $ is the initial ordinal, and 
$ k^+ $ is the smallest initial ordinal greater than $ k $. We note that an ordinal of the form $ k^+ $  is always regular ordinal.
In this paper we prove the following theorem.

\textbf{Тheorem 1.} {\it Let $\tau$ be an arbitrary initial regular ordinal, $\sigma$ and $\lambda$ be initial ordinals satisfying the inequality $\omega \leqslant \sigma <\lambda \leqslant \tau $. Then the space $C_p ([0, \tau \cdot \sigma])$ is not homeomorphic to the space $C_p ([0,\tau \cdot \lambda])$.}

If we combine this result with the results of \cite{Gorak},  we get a complete topological classification
of the spaces $C_p([0, \alpha])$ (which coincides with the uniform classification). We can write it in the form of the following theorem.

\textbf{Theorem 2.} {\it Let $\alpha$ and $\beta$ be ordinals and $\alpha\leq \beta$.

{\rm (a)} If $|\alpha|\neq |\beta|$, then 
	 $C_p ([0,\alpha])$ and $C_p([0,\beta])$ are not homeomorphic.

{\rm (b)}  If $\tau$  is an initial ordinal, $|\alpha|=|\beta|=\tau$  and either $\tau=\omega$ or $\tau$ is a singular ordinal or $\beta\geq\alpha\geq\tau^2$, then the spaces $C_p ([0,\alpha])$ and $C_p([0,\beta])$ are (uniformly) homeomorphic.

{\rm (c)} if $\tau$ is a regular uncountable ordinal and  $\alpha,\beta\in [\tau,\tau^2]$, then the space $C_p ([0,\alpha])$  is (uniformly) homeomorphic to the space $C_p([0,\beta])$  if and only if $\tau\cdot\sigma\leq\alpha\leq\beta<\tau\cdot\sigma^+$, where $\sigma$  is the initial ordinal, $\sigma<\tau$, and $\sigma^+$  is the smallest initial ordinal, exceeding $\sigma$.}

\begin{center}
\textbf{2. Proof of Theorem 1.}
\end{center}

We need some notation and auxiliary statements. For an arbitrary ordinal $\alpha$ and the initial ordinal $\lambda \leqslant \alpha$ we set
\[
A_{\lambda, \alpha} = \{t \in [0, \alpha] : \chi (t) = | \lambda| \},
\]
where
$\chi (t)$ is the character of the point $t \in [0, \alpha]$.
In particular, $A_{\omega, \alpha}$ is the set of all limit points of $t \in [0, \alpha]$,
having a countable base of neighborhoods.

Let $\alpha$ be a limit ordinal. The smallest order type of sets
$A \subset [0, \alpha]$ cofinal in $[0,\alpha)$, is called \textit{cofinality} of the ordinal $\alpha$ and denoted by $\mathrm{cf} (\alpha)$.

It is easy to see that $|\mathrm{cf} (\alpha)| = \chi (\alpha)$ for the limit ordinal $\alpha$.
The initial ordinal $\alpha$ is called \textit{regular} if $\mathrm{cf} (\alpha) = \alpha$. Otherwise, the initial ordinal is called \textit{singular}.

The symbol $D(x)$ denotes the set of points of discontinuity of the function $x$.

The proof of the following two lemmas is standard (see Example 3.1.27 in  \cite{Engel}).

\textbf{Lemma 1.} Let $\alpha$ be an arbitrary ordinal and let  $\tau$ be an initial ordinal such that $\omega <\tau\leqslant \alpha$, $t_0 \in A_{\tau, \alpha}$ and a function $x \colon [0, \alpha] \to \mathbb{R}$ is continuous at all points of the set $ A_{\omega, \alpha}$. Then there is an ordinal $\gamma <t_0$ such that $ x|_{(\gamma, t_0)} = \text{const}$. $\Box$

\textbf{Lemma 2.} If a function $x \colon [0, \alpha] \to \mathbb{R}$ is continuous at all points of the set $A_{\omega, \alpha}$, then the set $D(x)$ is at most countable. $\Box$

For the function $x\in \mathbb{R}^{[0, \alpha]}$ and the initial ordinal $\lambda \leqslant \alpha$ the symbol $G_\lambda (x)$ denotes the family
\[
G_\lambda (x) = \left\{ \bigcap_{s\in S} V_s : V_s\text{ is standard nbhd of   } x \text{ in } \mathbb{R}^{[0, \alpha]} \text{ and } |S| = |\lambda| \right\}.
\]
The elements of the family $G_\lambda (x)$ will be called \textit{$\lambda$ - neighborhoods of the function $x$}.

For a regular ordinal $\tau \geqslant \omega_1$ and a initial ordinal $\sigma \leqslant \tau$
we put
\begin{multline*}
M_{\tau \sigma} = \left \lbrace x \in \mathbb{R}^{[0, \tau \cdot \sigma]} : x \text {   is continuous at those points } t \in [0, \tau \cdot \sigma], \right. \\
\left. \vphantom {\mathbb{R}^{[0, \tau \cdot \sigma]}} \text{ for which } \mathrm{cf}(t) <\tau \right \rbrace.
\end{multline*}

It is clear, that $C([0, \tau \cdot \sigma]) \subset M_ {\tau \sigma}$.

\textbf{Lemma 3.} Let $\tau \geqslant \omega_1$ be an initial regular ordinal and let $\sigma$ be an initial ordinal such that  $\sigma \leqslant \tau$. Then 
\begin{multline*} 
M_{\tau \sigma} = \left\lbrace x\in \mathbb{R}^{[0,\tau \cdot \sigma]} : V \cap C_p([0,\tau \cdot \sigma]) \ne\varnothing  \text{ for every } V\in G_\lambda (x) \right. 
\\\left. \vphantom {\mathbb{R}^{[0,\tau\cdot\sigma]}} \text{ and each }\lambda < \tau \right\rbrace.
\end{multline*}

\begin{proof} We denote by $L_{\tau\sigma}$ the right-hand side of the equality and assume that
$x\notin M_{\tau\sigma}$, that is, $ x $ is discontinuous at some point $ t_0 $ for which $\mathrm{cf}(t_0) < \tau$. Since $ | \mathrm {cf}(t_0)| = \chi (t_0)$, there exists a base  $\{U_j(t_0)\}_{j\in J}$ of neighborhoods of the point $t_0$ such that $ | J | <\tau $. Since $ x $ is discontinuous at $ t_0 $, there exists a number $ \varepsilon_0> 0 $ such that for each
 $ j \in J $ there is a point $ t_j \in U_j (t_0)$ such that $ | x (t_j) -x (t_0) | \geqslant \varepsilon_0 $.
Let $ V = \bigcap \{V (x, t_j, t_0,1 / n) : j \in J, n \in \mathbb{N} \}$, where $ V (x, t_j, t_0,1 / n)$ is the standard neighborhood of the function $ x $ in the space $ \mathbb{R}^{[0, \tau \cdot \sigma]}$. If $ y \in V$, then $ y (t_j) = x (t_j) $ and $ y (t_0) = x (t_0) $.
Hence, the function $ y $ is discontinuous at the point $ t_0 $ and then
$ y \notin C_p ([0, \tau \cdot \sigma]) $. Thus, $ V \cap C_p([0, \tau \cdot \sigma]) = \varnothing $, that is, $ x \notin L_{\tau \sigma} $.

Now let $x \in M_{\tau \sigma}$, i.e. the function $x$ can be discontinuous only at the points of the set $A_{\tau, \tau \cdot \sigma}$. It is easy to see that the set
$ A_{\tau, \tau \cdot \sigma}$ has the form
 \[
A_{\tau,\tau\cdot\sigma}=\{\tau\cdot (\xi+1) : 0\leqslant\xi <\sigma\},  \text{ or }
\]
\[
A_{\tau,\tau\cdot\sigma}=\{\tau\cdot (\xi+1) : 0\leqslant\xi <\tau\}\cup \{\tau\cdot\tau\}, \text{ if } \sigma =\tau. 
\]

By Lemma 2, the set $D(x)$ is at most countable and therefore
\[
A_{\tau,\tau\cdot\sigma}\cap D(x)=
\{\tau \cdot(\xi_n+1) : \xi_n < \sigma, n\in\mathbb{N}\},  \text{ or } 
\]
\[
A_{\tau,\tau\cdot\sigma}\cap D(x)=
\{\tau\cdot(\xi_n+1) : \xi_n <\tau, n\in\mathbb{N}\}\cup \{\tau\cdot\tau\},  \text{ if } \sigma =\tau.
\]

Let $\lambda <\tau$ and  $V(x) = \bigcap\{U(x,\eta,1/n) : \eta\in S,\,n\in\mathbb{N}\}$ be a $\lambda$ -neighbourhood of the point $x$. Then $|S| <|\tau|$.
	
Since the countable set $A_{\tau,\tau\cdot\sigma}\cap D(x)$  is not cofinal in the regular ordinal $\tau\geqslant\omega_1$, for each $n\in\mathbb{N}$  there is an ordinal $\gamma_n$ such that $\tau \xi_n <\gamma_n <\tau (\xi_n + 1)$ and
$(\gamma_n, \tau (\xi_n + 1)) \cap S = \varnothing$. In the case $\sigma = \tau$ there is also an ordinal $\gamma_0 <\tau^2$,  such that $(\gamma_0, \tau^2) \cap S = \varnothing$ and 
$(\gamma_0, \tau^2) \cap \{\tau (\xi_n + 1) \}_{n = 1}^\infty = \varnothing$.

Consider the function
\[
\tilde{x}(t)=\left\lbrace 
\begin{array}{ll}
x(\tau(\xi_n+1)), & \text{ if } t\in (\gamma_n, \tau(\xi_n+1));\\
x(\tau^2), & \text{ if } t\in (\gamma_0, \tau^2);\\ x(t), & \text{ otherwise }.
\end{array}\right. 
\]
It is not difficult to see that the function $\tilde{x}$ is continuous at all points $t\in [0, \tau\cdot\sigma]$, and since $\tilde{x}|_S = x|_S$, $\tilde{x}\in V(x)$, that is, 
$V(x)\cap C_p ([0,\tau\cdot\sigma])\neq\varnothing$ and therefore $x\in L_{\tau\sigma}$.
\end{proof}

If $X$ is a Tikhonoff space, then the symbol $\nu X$ denotes the Hewitt completion of the space $X$. The proof of the following lemma can be found in \cite{Engel}, p. 218.

\textbf{Lemma 4.} If $\varphi\colon X\to Y$ is a homeomorphism of Tikhonoff spaces, then there exists a homeomorphism $\tilde{\varphi}\colon \nu X\to\nu Y$ such that $\tilde{\varphi}(x)=\varphi (x)$ for each $x \in X$.

\textbf{Lemma 5.} Let $\alpha$ be an arbitrary ordinal. Then
\[
\nu(C_p([0,\alpha])) = \left\lbrace x\in\mathbb{R}^{[0,\alpha]}: x \text{ is continuous at all points of the set } A_{\omega,\alpha}\right\rbrace.
\]
\begin{proof} It is known (\cite{Tkachuk}, p. 382) that for an arbitrary Tikhonov space $X$ the space $\nu (C_p (X))$ coincides with the set of all strictly $\aleph_0$ -continuous functions from $X$ to $\mathbb{R}$. In this case, the function $f \in \mathbb{R}^X$ is called strictly $\aleph_0$ -continuous (\cite{Arch}),  if for any countable set $A\subset X$ there is a continuous function $g \in \mathbb{R}^X$ such that $f|_A = g|_A$.

Since for each countable set $A \subset [0,\alpha]$, its closure $\bar{A}$ is also countable, by the Tietze-Uryson theorem we obtain that the set of all strictly $\aleph_0$-continuous functions in $[0,\alpha]$ in $\mathbb{R}$ coincides with the set of all those functions that are continuous on each countable subset $A\subset [0,\alpha]$. It is easy to see that these are precisely all those functions that are continuous at all points of the set $A_{\omega, \alpha}$.
\end{proof}

\textbf{Corollary 6.} If $\tau\geqslant\omega_1$ is the initial regular ordinal and $\sigma\leqslant\tau$ is the initial ordinal, then $M_{\tau\sigma}\subset \nu(C_p([0,\tau\cdot\sigma]))$.

For the initial ordinal $\sigma$ we denote by $\Gamma_\sigma$ the discrete space of cardinality $|\sigma|$ and consider the space
\[
c_0(\Gamma_\sigma)=
\left\lbrace x\in\mathbb{R}^{\Gamma_\sigma} : \{t\in\Gamma_\sigma\colon |x(t)|\geqslant\varepsilon\}\text{ is finite for any } \varepsilon >0\right\rbrace.
\]

\textbf{Lemma 7.} Let $\tau\geqslant\omega_1$ be an initial regular ordinal, $\sigma\leqslant\tau$  be an initial ordinal. Then there exists a homeomorphic embedding $f\colon c_0(\Gamma_\sigma)\to M_{\tau\sigma}$ such that $f(0)=0$ and $f(x)\in M_{\tau\sigma}\setminus C_p ([0,\tau\cdot\sigma])$, if $x\ne 0$.
\begin{proof}
We enumerate the points of the set $\Gamma_\sigma$ by the ordinals $\xi\in [0,\sigma)$.Then $\Gamma_\sigma = \{t_\xi\}_{\xi\in [0,\sigma)}$. For each characteristic function $\chi_{\{t_\xi\}}\in c_0(\Gamma_\sigma)$ we put $f(\chi_{\{t_\xi\}}) = \chi_{\{\tau(\xi+1)\}}$.It is obvious that $\chi_{\{\tau(\xi+1)\}}\in M_{\tau\sigma}\setminus C_p ([0,\tau\cdot\sigma])$. It remains to extend the map $f$ in the standard way to the space $c_0(\Gamma_\sigma)$.
\end{proof}

\textbf{Lemma 8.} Let $\tau\geqslant\omega_1$ be an initial regular ordinal, $\sigma,\lambda$ be an initial ordinals and $\omega\leqslant\lambda <\sigma\leqslant\tau$. If $f\colon c_0(\Gamma_\sigma)\to M_{\tau\lambda}$ is an injective mapping such that  $f(0)=0$ and $f(x)\in M_{\tau\lambda}\setminus C_p ([0,\tau\cdot\lambda])$ for $x\ne 0$, then the map $ f $ is not continuous.
\begin{proof} Suppose that there exists a continuous map $f\colon c_0 (\Gamma_\sigma) \to M_{\tau\lambda}$ with the above-mentioned properties. As in Lemma 7, let $\Gamma_\sigma = \{t_\xi \}_{\xi \in [1,\sigma)}$. Since the space $c_0 (\Gamma_\sigma)$ is considered in the topology of pointwise convergence, any sequence of the form $\chi_{\{t_{\xi_n} \}}$ converges to zero in this space. Consequently, at each point $\gamma \in [0, \tau\cdot\lambda]$ only a countable number of functions $f (\chi_{\{t_\xi \}})$ is nonzero. Since by the condition 
$f(\chi_{\{t_\xi\}})\in M_{\tau\lambda}\setminus C_p ([0,\tau\cdot\lambda])$, each function $f(\chi_{\{t_\xi \}})$ is discontinuous at some point of the set $A_{\tau,\tau\lambda} \subset [0, \tau\cdot\lambda]$.

Let
\[
B_\gamma = 
\left\lbrace f(\chi_{\{t_\xi\}}) : f(\chi_{\{t_\xi\}})\text{ is discontinuous at a point } \tau(\gamma+1)\in A_{\tau,\tau\lambda}\right\rbrace. 
\]
Since $\bigcup_{\gamma <\lambda}B_\gamma = f\left( \{\chi_{\{t_\xi\}} : \xi <\sigma\}\right)$ and $|\lambda|=|A_{\tau,\tau\lambda}|<|\sigma|$, there is a point $\gamma_0 <\lambda$, such that  $|B_{\gamma_0}|=|\sigma|$. Since at the point $\tau(\gamma_0+1)$ only a countable number of functions from $B_{\gamma_0}$ are nonzero, without loss of generality we can assume that all functions from $B_{\gamma_0}$ at the point $\tau(\gamma_0+1)$ are equal to zero. By Lemma 1, for each function $f(\chi_{\{t_\xi\}})\in B_{\gamma_0}$ there exists an ordinal  $\gamma_\xi <\tau(\gamma_0+1)$ such that  $f(\chi_{\{t_\xi\}})|_{[\gamma_\xi, \tau(\gamma_\xi +1))}=\text{const}=C_\xi$. Since $|B_{\gamma_0}| = |\sigma|>\omega$, in $B_{\gamma_0}$ there is an uncountable family of functions for which $|C_\xi|\geqslant\varepsilon_0$. Consider the sequence $\{f(\chi_{\{t_{\xi_n}\}})\}_{n=1}^\infty$ of such functions and put $\gamma_0 = \sup \{\gamma_{\xi_n}: n=1,2,\ldots\}$. Since $\mathrm{cf}(\tau(\gamma_0+1)) >\omega$, $\gamma_0 <\tau(\gamma_0+1)$  and therefore $|f(\chi_{\{t_{\xi_n}\}})(t)|\geqslant \varepsilon_0$ for each $t\in (\gamma_0, \tau(\gamma_0+1))$. But this contradicts the fact that the sequence $\{f(\chi_{\{t_{\xi_n}\}})\}_{n=1}^\infty$ converges pointwise to zero.
\end{proof}

\begin{proof}[Proof of Theorem 1.]
Suppose that there exists a homeomorphism $\varphi\colon C_p([0,\tau\cdot\sigma])\to C_p([0,\tau\cdot\lambda])$. We can assume that $\varphi (0)=0$. By Lemma 4, there exists a homeomorphism $\tilde{\varphi}\colon \nu(C_p([0,\tau\cdot\sigma]))\to \nu(C_p([0,\tau\cdot\lambda]))$ such that $\tilde{\varphi}(C_p([0,\tau\cdot\sigma]))=C_p([0,\tau\cdot\lambda])$. By Corollary 6  $M_{\tau\sigma}\subset \nu (C_p([0,\tau\cdot\sigma]))$, and by Lemma 3 $\tilde{\varphi}(M_{\tau\sigma}) = M_{\tau\lambda}$. 
By Lemma 7  the mapping $\tilde{\varphi}\cdot f: c_0 (\Gamma_\sigma)\to M_{\tau\lambda}$   is continuous, $(\tilde{\varphi}\cdot f)(0)=0$   and $(\tilde{\varphi}\cdot f)(M_{\tau\sigma})\subset M_{\tau\lambda}\setminus C_p ([0,\tau\cdot\lambda])$  for $x\neq 0$.
In this case, the map $\tilde{\varphi}|_{c_0 (\Gamma_\sigma)}$  is a homeomorphism of the space $c_0 (\Gamma_\sigma)\subset M_{\tau\sigma}$ onto the subspace $M_{\tau\lambda}$ such that $\tilde{\varphi}(0)=0$ and $\tilde{\varphi}(x)\subset M_{\tau\lambda}\setminus C_p ([0,\tau\cdot\lambda])$  for $x\ne 0$. But this is impossible by Lemma 8.
\end{proof}

The authors are grateful to the anonymous referee for helpful comments and suggestions to improve the manuscript.

\bigskip

CONTACT INFORMATION
\medskip

L.V.~Genze: \\
Faculty of Mathematics and Mekhaniks  \\
        Tomsk State University  \\
         Pr. Lenina 36  \\
         634050 Tomsk, Russia \\genze@math.tsu.ru 

\medskip

S.P.Gul'ko: \\
Faculty of Mathematics and Mekhaniks  \\
        Tomsk State University  \\
         Pr. Lenina 36  \\
         634050 Tomsk, Russia \\gulko@math.tsu.ru

\medskip

T.E.Khmyleva: \\
Faculty of Mathematics and Mekhaniks  \\
        Tomsk State University  \\
         Pr. Lenina 36  \\
         634050 Tomsk, Russia \\tex2150@yandex.ru

\end{document}